\documentclass[12pt,twoside]{amsart}
\input{includeNice3}
\usepackage{mabliautoref}
\usepackage{amssymb}
\usepackage{amscd}
\usepackage[abbrev,alphabetic]{amsrefs}
\usepackage{hyperref}
\usepackage[a4paper,margin=1in]{geometry}  
\usepackage{CJKutf8}
\usepackage{soul}

\usepackage{xypic}
\usepackage{tikz-cd}
\usepackage{caption}

\usepackage{url}
\RequirePackage{booktabs, multirow}
\RequirePackage{pgf}

\title[Counterexamples to the MMP for 1-foliations in characteristic $p$]
{Counterexamples to the MMP for 1-foliations in positive characteristic} 
\author[F. Bernasconi]{Fabio Bernasconi} 

\address{Departement Mathematik und Informatik, Universit{\"a}t Basel, Spiegelgasse 1, CH-4051 Basel, Switzerland.} 
\email{fabio.bernasconi@unibas.ch}

\begin{document}

\maketitle

\begin{abstract}
We show that many statements of the Minimal Model Program, including the cone theorem, the base point free theorem and the existence of Mori fibre spaces, fail for 1-foliated surface pairs $(X,\mathcal{F})$ with canonical singularities in characteristic $p>0$.
\end{abstract}

\tableofcontents

\section{Introduction}
The Minimal Model Program (MMP) is a conjecture concerning the birational geometry of algebraic varieties, with the goal of extending the theory of minimal models from algebraic surfaces to higher dimensions. A large part of this program has been established for varieties over fields of characteristic 0 \cite{KM98, BCHM10, LM22}, as well as for 3-dimensional varieties in positive and mixed characteristic \cite{Tan18, HX15, 7authors, TY20}.

In recent decades, birational geometry of foliations in characteristic 0 has emerged as a highly active research area in algebraic geometry. 
Starting from the work of Brunella, Mendes and McQuillan, the MMP for foliations on surfaces has been proven in \cite{McQ08, Bru15} and applied to solve the Green--Griffiths conjecture for surface with positive second Segre class \cite{McQ98}. 
More recently, significant progress has been made in extending the MMP to foliations on 3-folds, as discussed in \cite{CS21, SS22, CS20}. Additionally, the MMP for algebraically integrable foliations in every dimension has been established in \cite{CS23, CHLX}.
Furthermore, the cone theorem for canonical Gorenstein foliations of rank 1 in every dimension has been established by Bogomolov and McQuillan \cite{BM16}*{Corollary 4.2.1}, and it has recently been extended to the log canonical setting in \cite{CS23b}

In light of the aforementioned developments, it is natural, from the perspective of birational geometry, to inquire whether one can extend the MMP to the case of foliations in characteristic $p>0$. The purpose of this brief note is to provide a negative response by constructing several counterexamples to the classical MMP statements for 1-foliations on surfaces in positive characteristic.

We begin with the cone theorem. 
Tanaka established the cone theorem for $\mathbb{Q}$-factorial projective surfaces over a perfect field of characteristic $p > 0$ in \cite{Tan14}*{Theorem 3.13}. 
In \cite{Spi20}*{Theorem 6.3}, Spicer demonstrated the cone theorem for arbitrarily singular foliations of rank 1 on normal projective surfaces over $\mathbb{C}$. 
Our first result shows that the cone theorem does not hold for 1-foliations in characteristic $p > 0$.

\begin{theorem} \label{thm: cone_counterexample}
For every algebraically closed field $k$ of characteristic $p>0$, there exists a foliated surface pair $(X, \mathcal{F})$ such that
\begin{enumerate}
    \item $\mathcal{F}$ is a canonical 1-foliation of rank 1 on $X$ such that $-\KF$ is ample;
    \item the Mori cone $\NE(X)$ is not finitely generated nor locally polyhedral;
    \item $X$ does not contain any rational curves.
\end{enumerate}
\end{theorem}

Another pillar of the MMP in characteristic 0 is the base point free theorem \cite{KM98}*{Theorem 3.3}.
For $\mathbb{Q}$-factorial surfaces in characteristic $p>0$, this has been proved in \cite{Tan14}*{Theorem 3.21}.
A version of the basepoint free theorem has been established in the case of rank 1 foliations on surfaces in characteristic 0 by McQuillan (\cite{McQ08}*{Section III.2}) and some special cases in higher dimensions are proved (\cite{CS21}*{Theorem 9.4} and \cite{CS20}*{Theorem 9.2}).
Here, we demonstrate that the base point free theorem does not generally hold for 1-foliations in characteristic $p > 0$.

\begin{theorem} \label{thm: bpf_fails}
    For every uncountable algebraically closed field $k$ of characteristic $p>0$, there exists a $\mathbb{Q}$-factorial canonical 1-foliated surface pair $(X, \mathcal{F})$ such that
    \begin{enumerate}
        \item $D$ is a nef $\mathbb{Q}$-divisor on $X$ such that $D-\KF$ is ample, and
        \item $D$ is not semi-ample. 
    \end{enumerate}
\end{theorem}

Finally, we demonstrate that Mori fiber spaces may not exist for foliated pairs with a non-pseudo-effective canonical divisor.

\begin{theorem} \label{thm: MFS_counterexample}
    For every algebraically closed field of characteristic $p>0$, there exists a canonical 1-foliated surface pair $(X, \mathcal{F})$ such that 
    \begin{enumerate}
        \item $\KF$ is not pseudo-effective, and
        \item it does not exists a $\KF$-negative birational contraction $X \to Y$ such that $Y$ admits a $\KF$-Mori fibre space structure.
    \end{enumerate}
\end{theorem}

The main source of examples comes from the natural foliation $\mathcal{F}$ associated to a $p$-cyclic covering $X=Y[\sqrt[p]{f}] \to Y$ in characteristic $p>0$, which we recall in \autoref{ss-pcyclic}. 
This construction, introduced by Ekedahl \cite{Eke87}*{page 145} and revisited by Koll\'{a}r in \cite{Kol95}, gives a counterexample to the Bogomolov--Sommese vanishing theorem in positive characteristic. It has recently been revisited by Langer to give counterexamples to Miyaoka's generic semipositivity \cite{Lan21}*{Section 5.1} and by Graf to construct counterexamples to the logarithmic extension theorem for differential forms (\cite{Gra21}*{Theorem 1.6}).
Applying this construction to $p$-cyclic covering of well-chosen abelian surfaces, we construct Fano foliations which serve as the sources of counterexamples for our results in \autoref{ss: counter}.

\begin{remark}
    The results of these notes suggest that general 1-foliations exhibit pathological behavior concerning the MMP. 
    However, in the case of foliations defined as infinitesimal groupoids (as per \cite{McQ22}*{Definition 1.1}), McQuillan has recently established a cone theorem for absolutely $\mathbb{Q}$-Gorenstein foliations of rank 1 (\cite{McQ22}*{Fact 5.7}) in any characteristic. This shows that \autoref{thm: cone_counterexample} represents a phenomenon unique to such pathological 1-foliations.
   We also observe that our examples are 1-foliations that are not $(\infty)$-foliations in the sense of \cite{Gab23} as explained in \cite{Eke87}*{page 145}. This leads us to inquire whether the MMP might still hold true in the case of $(\infty)$-foliations on surfaces.
\end{remark}

\textbf{Acknowledgements.}
I would like to thank P.~Cascini, A.~Langer, C.~Spicer, R.~Svaldi and P.~Grabowski  for interesting comments on the content of this note and J.V. Pereira for useful discussion on singularities and algebraicity of $1$-foliations.
The author was partly supported by the grants $\#200021/169639$ and PZ00P2-216108 from the Swiss National Science Foundation.

\section{Preliminaries}

\subsection{Notations}

\begin{enumerate}
    \item Throughout this article, we fix $p$ to be a prime number.
    \item We denote by $k$ an algebraically closed field.
    \item A $k$-variety (or simply variety) is a quasi-projective integral scheme of finite type over $k$.  
    \item If $X$ is a variety over a field $k$ of characteristic $p>0$, we denote by $F \colon X \to X$ the \emph{absolute Frobenius} of $X$.
    \item If $X$ is a projective normal $k$-variety, we denote by $\NS(X)$ the real numerical N\'eron-Severi space and by $\N1(X)$ the real vector space of numerical equivalence classes of curves.
    The \emph{Mori cone} $\NE(X)$ of $X$ is the closure of the real cone spanned by effective curves inside $\N1(X)$. 
    \item Given a normal $k$-variety $X$, we define its \emph{tangent sheaf} $T_X$ as the dual of $\Omega_{X/k}^{1}$ and we denote  $[ - , - ] \colon T_X \times T_X \to T_X$ the Lie bracket on $T_X$.
\end{enumerate}

\subsection{Foliations in characteristic $p>0$}

We recall the definition of a foliation on algebraic varieties.
 
\begin{definition}
	Let $X$ be a normal variety over $k$. 
    A \emph{foliation} $\mathcal{F}$ is a coherent subsheaf of $T_X$ such that
    \begin{enumerate}
        \item it is saturated in $T_X$, i.e. $T_X/\mathcal{F}$ is torsion-free;
        \item it is closed under Lie brackets, i.e. $[\mathcal{F}, \mathcal{F}] \subset \mathcal{F}$.
    \end{enumerate}
We say $\mathcal{F}$ is a foliation of rank $r$ if the rank of the coherent sheaf  $\mathcal{F}$  is $r$. 
A \emph{canonical divisor} $\KF$ of a foliation $\mathcal{F}$ is a Weil divisor such that $\mathcal{O}_X(\KF) \simeq \det(\mathcal{F})^{*}.$
 We define $\mathcal{N}_\mathcal{F}^* = (T_X / \mathcal{F})^*$ to be the \emph{conormal sheaf} of $\mathcal{F}$. 
 The normal sheaf $\mathcal{N}_\mathcal{F}$ of $\mathcal{F}$ is the dual of the conormal sheaf.
\end{definition}

From now on, we assume $k$ to be a field of characteristic $p>0$.
In this case, together with Lie brackets, we have an additional operation on vector fields of a variety $X$: the elevation to the $p$-th power. 
More precisely, we define a homomorphism of sheaves of abelian groups
\[\circ^{[p]} \colon F^*T_X \rightarrow T_X, \] 
as follows: given an open set $U \subset X$ and for every $\mathcal{O}_X(U)$-derivation $D \in T_X(U)$, its $p$-th power 
$D^{[p]}$ is the composition $\underbrace{D \circ \dots \circ D}_{p\text{ times}}$.
Note that one can verify that $D^{[p]}$ is still a derivation by applying the Leibniz formula and the fact that $p=0$.
It is easy to see that, for any $f \in \mathcal{O}_X$ and $D \in T_X$, then $(fD)^{[p]}=f^pD^{[p]}$. 
Note that the elevation to the $p$-th power does not respect the additive structure.
In this note, we work with foliations invariant under the elevation to the $p$-th power.

\begin{definition}
	Let $X$ be a normal variety over $k$. A foliation $\mathcal{F}$ on $X$ is a \emph{1-foliation} (or a $p$-closed foliation) if the composite map
	$$F \colon F^*\mathcal{F} \rightarrow F^*T_X \xrightarrow{\circ^{[p]}} T_X \rightarrow TX/\mathcal{F}$$
	vanishes identically.
    Equivalently, for every open set $U \subset X$ and $D \in \mathcal{F}(U)$, the $p$-th power $D^{[p]} \in \mathcal{F}(U)$.
\end{definition}

\begin{remark}
    Note that the composition $F^*\mathcal{F} \rightarrow F^*T_X \rightarrow T_X \rightarrow TX/\mathcal{F}$ is $\mathcal{O}_X$-linear.
\end{remark}

1-foliations appear naturally in the Jacobson correspondence for purely inseparable field extension of height 1 first proven in \cite{Jac44} (and later extended in \cite{Jac89}). 
For the formulation in terms of factorisation of the geometric Frobenius we refer to 
\cite{PW22}*{Proposition 2.9} (see also \cite{Eke87}).

\begin{theorem}[Jacobson's correspondence] \label{thm: jacobson_corr}
Let $X$ be a normal variety over $k$. 
There is a one-to-one correspondence between
    \begin{enumerate}
        \item 1-foliations $\mathcal{F}$ of rank r; 
        \item Factorisation of  the geometric Frobenius  morphism $F_{X/k} \colon X \xrightarrow{f} Y \xrightarrow{g} X^{(-1)}$.
    \end{enumerate}
If $\mathcal{F}$ has rank r, then the degree of $f$ is $\deg(f)=p^r$.
\end{theorem}

We recall the formula for the ramification formula for the canonical classes in the case of purely inseparable morphism induced by a $1$-foliation.
 
\begin{proposition}[{\cite{PW22}*{Proposition 2.10}}] \label{ramificationquotfoliations}
	Let $ \mathcal{F}$ be a 1-foliation on a normal variety $X$. 
    Let $ \pi \colon X \rightarrow X/\mathcal{F}$ be the purely inseparable morphism induced by $\mathcal{F}$ under the Jacobson correspondence. Then
	\[ \pi^* K_{X/\mathcal{F}} - K_X = (p-1) K_{\mathcal{F}}. \] 
\end{proposition}

An instance where 1-foliations of rank $r$ appear naturally is in presence of a separable dominant morphism $f \colon X \to Y$ whose generic fibre has dimension $r$. In this case, the foliation $\mathcal{F} = \text{ker}(df \colon T_X \rightarrow f^*T_Y)$ is a 1-foliation of rank $r$.

\begin{definition}\label{def: int_foliations}
	A foliation $\mathcal{F}$ on $X$ is said to be \emph{algebraically integrable} if there exists an open Zariski subset $U \subset X$ and a smooth morphism $f \colon U \rightarrow V$ such that $\mathcal{F}|_U = \text{ker}(df \colon TU \rightarrow f^*TV)$. 
\end{definition}

\begin{remark}
As explained in \cite{LPT18}*{Section 7, pages 1362-1363} the class of $1$-foliations coincides with the class of algebraically integrable foliations in characteristic $p>0$.
\end{remark}

For smooth 1-foliations on smooth varieties, we have an analogue of Frobenius theorem.

\begin{theorem}[\cite{McQ08}*{Divertimento II.1.6.}] \label{thm: Frobenius}
    Let $\mathcal{F}$ be a smooth 1-foliation of rank 1 on a smooth variety $X$ of dimension $n$ (eq. $\mathcal{F} \to T_X$ is a morphism of vector bundles). 
    Then at each point $x$, there is a choice of local parameters $x_1, \dots x_n$ around $x$ such that the vector field $\partial_{x_1}$ generates $\mathcal{F}$.
\end{theorem}


Given a foliation $ \mathcal{F}$ on $X$ and a proper birational morphism $\pi \colon Y \to X$, let $\widetilde{\mathcal{F}}$ be the pull-back foliation on $Y$ defined as in \cite{Dru18}*{Section 3.3}. 
In the spirit of the MMP, one can introduce the class of singularities for foliated pairs by comparing $\KF$ and $K_{\widetilde{\mathcal{F}}}$.
For this, it is natural to require the  canonical class of the foliation to be $\mathbb{Q}$-Cartier.

\begin{definition}\label{def: foliated_pair}
    We say $(X, \mathcal{F})$ is a \emph{foliated pair} (resp. \emph{1-foliated pair}) if $X$ is a normal variety, $\mathcal{F}$ is a foliation (resp. 1-foliation) such that the canonical divisor $\KF$ is $\mathbb{Q}$-Cartier. 
\end{definition}

As usual, we start with the definition of discrepancies.
If $(X, \mathcal{F})$ is a foliated pair and $\pi \colon Y \to X$ is a proper birational morphism of normal varieties, We can write
$$K_{\widetilde{\mathcal{F}}} = \pi^*\KF +\sum_i a(E_i, \mathcal{F}) E_i, $$
where $E_i$ runs through the $\pi$-exceptional divisors. 
We say that $a(E, \mathcal{F})$ is the \emph{(foliated) discrepancy} of $E$ with respect to $\mathcal{F}$.

\begin{definition}
We say a prime divisor $D \subset X$ is an \emph{invariant} divisor for the foliation $\mathcal{F}$ if  $\mathcal{F}|_D \to T_X|_D$ factorises through $T_D$ at the generic point of $D$. 
For each prime divisor $D$ on $X$, we define $\epsilon(D)=0$ (resp. $\epsilon(D)=1$)  if $D$ is invariant (resp. not invariant) for $\mathcal{F}$.
\end{definition}

\begin{definition}
    A foliated pair $(X, \mathcal{F})$ has \emph{terminal} (resp. \emph{canonical}) singularities if $a(E, \mathcal{F}) >0 $ (resp. $\geq 0$) for every exceptional prime divisor $E$ over $X$.
    It has \emph{klt} (resp. \emph{log canonical}) if $a(E, \mathcal{F}) >-\epsilon(E) $ (resp. $\geq -\epsilon(E)$) for every exceptional prime divisor $E$ over $X$.
\end{definition}

\begin{remark}\label{rem: can_sing}
    Let $\mathcal{F}$ be a Gorenstein foliation on $X$ generated by a vector field $v$. 
    By definition, $\mathcal{F}$ is canonical if and only if for every proper birational morphism $f \colon Y \to X$ the rational vector field $\widetilde{v}$ obtained as pull-back of $v$ does not vanish at the generic point of any exceptional divisor.
\end{remark} 

The next two examples should warn the reader that the notions of terminal and canonical singularities are more subtle in characteristic $p>0$ than in characteristic 0.

\begin{example} \label{ex: smooth_fol}
    Smooth 1-foliations $\mathcal{F}$ on smooth surfaces are canonical.  
    By \autoref{thm: Frobenius} we can suppose $\mathcal{F}$ is the smooth foliation generated by the vector field $\partial_x$ in $\mathbb{A}^2_{x,y}$.
    For every proper birational morphism $f \colon Y \to X$, let $v$ be the rational vector field obtained as the pull-back of $\partial_x$.
    As $f$ is birational, we have that $f^*(dx)(v)=1$. As $f^*(dx)$ is a regular 1-form, we conclude that $v$ cannot acquire zeroes at the generic point of any exceptional divisors, showing that $\mathcal{F}$ has canonical singularities by \autoref{rem: can_sing}.
    
	However, it is important to note that $\mathcal{F}$ is not terminal in characteristic $p > 0$, as pointed out in \cite{McQ08}*{Counterexample I.2.11}. This contrasts with the case of terminal excellent surfaces, which are shown to be regular in \cite{kk-singbook}*{Theorem 2.29}. 
    To show that $(X, \mathcal{F})$ is not terminal, it is sufficient to perform the weighted blow-up $\pi \colon Y \coloneqq \text{Bl}_{(1,p)}\mathbb{A}^2 \rightarrow \mathbb{A}^2$ with weight $(1, p)$ at the origin. 
    In a local chart of $Y$ containing the generic point of the exceptional divisor $E$ of $\pi$, the morphism $\pi$ is described by $(s,t) \mapsto (s, s^p t).$
    As $dy$ is a local section of the conormal bundle $\mathcal{N}^*_{\mathbb{A}^2}$
    together with the formulas 
	$\pi ^*dy=s^pdt, \text{ and }\pi^*(dx \wedge dy)=s^p ds \wedge dt$
	we conclude that
    $$K_{Y}  = \pi^{*}K_X+pE \text{ and } \mathcal{N}^*_{Y} = \pi^* \mathcal{N}^*_{\mathbb{A}^2}+pE.$$ 
    Thus $K_{\widetilde{\mathcal{F}}} = \pi^* K_{\mathcal{F}}$, concluding that $\mathcal{F}$ does not have terminal singularities. 
\end{example}

If $(X,\mathcal{F})$ is a canonical surface foliated pair in characteristic 0, then $X$ has log canonical singularities \cite{McQ08}*{Fact I.2.4 and Theorem III.3.2}, and this has been extended to rank 1 foliations on 3-folds in \cite{CS21}*{Theorem 1.5}. 
The following shows that this is no longer true for canonical 1-foliations in positive characteristic.

\begin{example} \label{ex: A_{p-1}-foliations}
    Let $X:= \left\{ z^p-f(x,y)=0 \right\} \subset \mathbb{A}^3_{x,y,z}$ be a normal hypersurface singularity. 
    In characteristic $p>0$, $X$ is an invariant divisor on $\mathbb{A}^3$ for the 1-foliation described by the vector field $\partial_z$ on $\mathbb{A}^3$.  
    Therefore we can restrict $\partial_z$ on $X$, which naturally defines an induced 1-foliation $\mathcal{F}$ on $X$. Note that in characteristic 0, $X$ is not invariant for $\partial_z$.
    
    We claim that $\mathcal{F}$ has canonical singularities. 
    Let $f \colon Y \to X$ be a proper birational morphism and let $v$ be the rational vector field on $Y$ obtained as the pull-back $\partial_z$.
    As $f^*(dz)$ is a regular 1-form on $Y$ and $f^*(dz)(v)=1$, we conclude that $v$ does not acquire zeroes along the exceptional divisor and $\mathcal{F}$ has canonical singularities by \autoref{rem: can_sing}. 
    Note that the singularities of $X$ can be easily worse than log canonical by choosing appropriate $f$.
\end{example}

\section{Proof of the main theorems}

In this section, we show the main counterexamples to the MMP for 1-foliations of this note.
We start with the review the construction of $p$-cyclic covering and their properties in characteristic $p>0$ following \cite{Kol95}.

\subsection{1-foliations from $p$-cyclic covering} \label{ss-pcyclic}

We fix $k$ to be an algebraically closed field of characteristic $p>0$ and let $Y$ be a normal projective variety over $k$ of dimension $n$ endowed with an invertible sheaf $L$. 

Let $\pi \colon \mathbb{L}:=\Spec_Y \bigoplus_{m \geq 0} L^{-\otimes m} \rightarrow Y$ be the associated line bundle to $L$ together with its tautological section $y_L \in H^0(\mathbb{L}, \pi^*L)$.
Given $s \in H^0(Y, L^{\otimes d})$ a non-zero global section for some $d>0$, we can construct the $d$-th cyclic cover $\pi \colon X \rightarrow Y$ along $s$ defined by the vanishing of the global section $(y_L^{d}-\pi^*s) \in H^0(\mathbb{L}, \pi^*L^{\otimes d})$. Sometimes we denote $X$ by $Y[\sqrt[d]{s}]$. 
If $U \subset Y$ is an open subset such that $L|_{U} \simeq \mathcal{O}_U$, we have an isomorphism $X \simeq \Spec_U \mathcal{O}_U[y]/(y^d-f)$ over $U$ for some $f \in \mathcal{O}_U$. 
In particular, if $Y$ is Cohen--Macaulay (resp. Gorenstein, locally complete intersection), then so is $X$. 
We recall a few basic sequences of differentials we will need to construct a 1-foliation on $X$. 

\begin{lemma}[cf. \cite{Kol95}*{Lemma 9}] \label{ovvio}
    We have the following short exact sequences:
    \begin{enumerate}
        \item On $\mathbb{L}$, we have 
	$0 \rightarrow \pi^*\Omega_Y^1 \rightarrow \Omega^1_\mathbb{L} \rightarrow \pi^*L^{-1} \rightarrow 0. $	
 \item  On $X$, we have $\pi^*L^{-\otimes d} \xrightarrow{d_X} \Omega^1_\mathbb{L}|_X \rightarrow \Omega^1_X \rightarrow 0,$ where, on the smooth locus of $X$, $d_X$ is locally described by
	\[ \left( -\frac{\partial s}{\partial x_1} dx_1, \dots, -\frac{\partial s}{\partial x_n} dx_n, dy^{d-1}dy  \right). \] 
    \end{enumerate}
\end{lemma}

We now specialize to the case where $p$ divides $d$. 

\begin{lemma}\label{l-fol-constructed}
    If $p \mid d$, we have the following exact sequence on $X$:
	\[ 0 \rightarrow \coker \left(\pi^*L^{-\otimes d} \rightarrow \pi^*\Omega^1_Y \right) \rightarrow \Omega_X^1 \rightarrow \pi^*{L^{-1}} \rightarrow 0. \]
\end{lemma}


Let $Y$ be a smooth variety together with a very ample line bundle $L$ and a global section $s \in H^0(Y, L^{\otimes p})$. 
On $X \coloneqq Y[\sqrt[p]{s}] $, the homomorphism $\Omega_X^1 \to \pi^*L^{-1}$ constructed in \autoref{l-fol-constructed} defines a foliation $\mathcal{F}$ of rank 1 on $X$ whose canonical bundle $K_{\mathcal{F}}=\pi^*{L}^{-1}$ is anti-ample. 
We say $\mathcal{F}$ is the \emph{foliation associated to the $p$-th cyclic cover} defined by the pair $(L, s)$.

\begin{lemma} \label{lem: p-closed-check}
    The foliation $\mathcal{F}$ is $p$-closed.
\end{lemma}

\begin{proof}
   In analytic coordinates, $X$ is described by $\left\{ z^p-f=0 \right\}$ where $f \in k[x_1, \dots, x_n]$.
   The foliation $\mathcal{F}$ is generated by the vector field $\partial_z$, which is obviously $p$-closed.    
\end{proof}

In general, $X$ is not smooth even if ($s=0$) is a smooth divisor. 
Nevertheless, the singularities appearing on $X$ can be controlled by choosing $s$ generically enough as explained in \cite{Kol95}.
When $p \mid d$, then $L^{\otimes d}$ has a natural operator:
\[ d \colon L^{\otimes d} \rightarrow L^{\otimes d} \otimes \Omega^1_Y,  \]
which is defined as follows (see \cite{Kol95}*{Definition-Lemma 13}): choose a local generator $\tau$ of  $L$ and for each section $s =f \tau^{\otimes d}$ of $L^{\otimes d}$ we define $d(s) \coloneqq df \otimes \tau^{\otimes d}$. This operator is indipendent of the choice of a generator $\tau$.
If we choose a local system of parameters $x_1, \dots, x_n$ near a point $x$, the matrix $$H(s) \coloneqq \left( \frac{\partial^2 f}{\partial x_i \partial x_j} \right)_{i, j}$$ is called the Hessian of $s$ and the rank of $H(s)(x)$ is independent of the choices made.

\begin{definition}
    We say that $s \in H^0(X, L^{\otimes d})$ has a \emph{critical point} at $x$ if $ds(x)=0$. 
    We say a critical point $x$ is \emph{non-degenerate} if the rank of Hessian $H(s)(x)$ is equal to the dimension of $X$.
\end{definition}

We recall the following local result which relates the zeros of $d s$ with the singularities of $Y[\sqrt[k]{s}]$ and describes the singularities of the cover if the points are non-degenerate.

\begin{lemma}[\cite{Kol95}*{20}] \label{lem: localcomp}
	If $p \mid d$, the hypersurface ($y^d-s(x_1, \dots, x_n)=0$) is singular at the point $(y, x_1, \dots, x_n)$ if and only if $(x_1, \dots, x_n)$ is a critical point of $s$.
    Moreover, $s$ is non-degenerate if and only if we can write, up to a change of coordinates, 
$$s(x_1, \dots, x_n)=
\begin{cases} x_1x_2+x_3x_4+ \dots +x_{n-1}x_n +f_3(x), & \mbox{if }n\mbox{ even} \\ 
x_1^2+x_2x_3+ \dots + x_{n-1}x_n+f_3(x), & \mbox{if }n\mbox{ odd} 
\end{cases}$$
where $f_3 \in (x_1, \dots, x_n)^3$.
\end{lemma}

\subsection{Counterexamples} \label{ss: counter}
 
We recall how the numerical spaces of cycles change under purely inseparable morphisms.

\begin{lemma} \label{lem: NSunderInseparable}
	Let $f \colon X \rightarrow Y$ be a purely inseparable morphism of projective normal varieties over $k$. Then $f_* \colon \N1(X) \rightarrow \N1(Y)$ is an isomorphism. Under this isomorphism, $\NE(X)$ is mapped onto $\NE(Y)$.
\end{lemma}

\begin{proof}
    The proof is dual to the one in \cite{Kee99}*{Lemma 1.4}.
    As $f$ is purely inseparable, there exists $l \geq 0$ and a factorisation of the $l$-th geometric Frobenius morphism $F^{l}_{X/k} \colon X \xrightarrow{f} Y \to X^{(-l)}$. 
    Given a curve class $[C] \subset \N1(X)$, we have $(F^l_{X/k})_* [C]=p^l [C]$ and this concludes the proof of the lemma.
\end{proof}

We are now ready to introduce the main example of this note by taking $p$-cyclic covering of abelian surfaces. A similar construction has been used in \cite{Lan21}*{Section 5.1} to construct counterexamples to Miyaoka's semipositivity theorems in positive characteristic.

\begin{proposition} \label{prop: Fano-foliations}
    Let $Y:=A$ be an abelian surface equipped with a very ample line bundle $L$ on $Y$, and let $s \in H^0(X,L^{\otimes p})$ be a general section. 
    Let $X:=Y[\sqrt[p]{s}]$ be the $p$-th cyclic covering together with the induced foliation $\mathcal{F}$.
    Then
    \begin{enumerate}
        \item $X$ has canonical singularities of type $A_{p-1}$,
        \item $(X, \mathcal{F})$ is a canonical 1-foliated pair,
        \item $-\KF$ is ample (i.e. $\mathcal{F}$ is a Fano foliation with canonical singularities).
    \end{enumerate}
\end{proposition}

\begin{proof}
    By \cite{Kol95}*{Proposition 18}, we can choose $s$ to have only non-degenerate critical points.
    The fact that $\mathcal{F}$ is a 1-foliation follows from \autoref{lem: p-closed-check} and the Fano condition follows from $\KF=\pi^{*}L^{-1}$.
	On a smooth big open set $U$ of $Y$ along which $s$ has no critical points, the foliation $\mathcal{F}$ is generated by the vector field $\partial_z$, thus concluding that $\mathcal{F}$ is a smooth foliation on $p^{-1}(U)$, and so canonical by \autoref{ex: smooth_fol}.
    Around a non-degenerate critical point, we have in analytic coordinates that $X$ is described by $\left\{z^p-xy=0 \right\}$ by \autoref{lem: localcomp} and the foliation $\mathcal{F}$ is thus canonical by \autoref{ex: A_{p-1}-foliations}.
\end{proof}

\begin{proof}[Proof of \autoref{thm: cone_counterexample}] \label{conethmfail}

Let $Y$ be an abelian surface of Picard rank $\rho(Y) \geq 3$, together with an ample line bundle $L$ and consider the foliated pair $(X, \mathcal{F})$, where $\mathcal{F}$ is the Fano 1-foliation of \autoref{prop: Fano-foliations}. 
Note there are no rational curves on $X$ since $X$ admits a finite morphism to an abelian surface. 
Moreover, the Mori cone $\NE(X)$ is isomorphic to $\NE(Y)$ by \autoref{lem: NSunderInseparable} and this is not finitely generated nor locally polyhedral by \cite{Deb01}*{6.3}.
\end{proof}

\begin{remark}
	A bend-and-break result  for $1$-closed foliation has been proven by Langer (\cite{Lan15}*{Theorem 2.1}). 
    In order to produce rational curves one needs to require the following lower bound on the foliated canonical class:
	\[-K_\mathcal{F} \cdot C > \frac{K_X \cdot C}{p-1}. \]
	Note that our example shows that the bound of Langer for the bend-and-break is optimal. 
    Indeed, the ramification formula \autoref{ramificationquotfoliations} 
	implies
	\[-K_{\mathcal{F}} =\pi^*L= \frac{K_X-\pi^*K_Y}{p-1}=\frac{K_X}{p-1}. \]
    Note that, even in Langer's range, we still do not know if the rational curves produced are tangent to the foliation.
\end{remark}

We show the base point free theorem fails for the case of 1-foliations.

\begin{proof}[Proof of \autoref{thm: bpf_fails}]
	Let $k$ be an uncountable algebraically closed field of characteristic $p>0$.
    Let $Y$ be an abelian surface over $k$, together with an ample line bundle $L$ and consider the foliated pair $(X, \mathcal{F})$, where $\mathcal{F}$ is the Fano 1-foliation of \autoref{prop: Fano-foliations}. 
	Let $D \in \Pic^0(A)$ be a nef divisor which is not semi-ample (equiv. not-torsion, whose existence is guaranteed by the hypothesis on $k$).
	By \cite{Kee99}*{Lemma 1.4}, on $X$ the divisor $\pi^*D$ is nef but not semi-ample. 
    Note that $\pi^*D - \KF \equiv -\KF$ is ample, thus concluding.
\end{proof}

We show that Mori fibre spaces do not always exist.

\begin{proof}[Proof of \autoref{thm: MFS_counterexample}]		
    Let $A$ be a simple abelian surface of Picard rank 2. By \cite{Bau98}*{Lemma 1.1 and Proposition 1.2} (the characteristic assumption is not needed in the proof), the two extremal rays of the nef cone are irrational.  
    Let $X:=A[\sqrt[p]{s}]$ be a $p$-cyclic cover with the foliation $\mathcal{F}$ of \autoref{prop: Fano-foliations}. 
    In this case, by \autoref{lem: NSunderInseparable} the two extremal rays of the nef cone of $X$ are irrational and thus $X$ does not admit any Mori fibre space structure.
\end{proof}

We conclude by noting the BAB conjecture fails for canonical Fano 1-foliations.

\begin{example}[Non-boundedness of Fano] \label{ex: BAB-fails}
	The foliated pairs $(X, \mathcal{F})$ obtained by taking larger and larger powers of $L$  in \autoref{prop: Fano-foliations} show that canonical Fano 1-foliations on surfaces in characteristic $p>0$ are not foliated bounded as there is no upper bound to the self-intersection numbers $\KF^2$.
\end{example}
 
\bibliographystyle{amsalpha}
\bibliography{refs}

 \end{document}